\documentclass[aps,prd,preprint,groupedaddress]{revtex4-1}

\usepackage{graphicx}
\usepackage{dcolumn}
\usepackage{bm}
\usepackage{amsmath}
\usepackage{epstopdf}

\begin{document}


\title{Modified Iterated Crank-Nicolson Method with Improved Accuracy}


\author{Qiqi Tran}
\email{qhuang11@students.desu.edu}
\author{Jinjie Liu}
\email{jliu@desu.edu}
\affiliation{Department of Mathematical Sciences, Delaware State University}


\date{\today}

\begin{abstract}
The iterated Crank-Nicolson (ICN) method is a successful numerical algorithm in numerical relativity for solving partial differential equations. 
The $\theta$-ICN method is the extension of the original ICN method where $\theta$ is the weight when averaging the predicted and corrected values. 
It has better stability when $\theta$ is chosen to be larger than 0.5, but the accuracy is reduced 
since the $\theta$-ICN method is second order accurate only when $\theta$ = 0.5. 
In this paper, we propose two modified $\theta$-ICN algorithms that have second order of convergence rate when $\theta$ is not 0.5, 
based on two different ways to choose the weight $\theta$. 
The first approach employs two geometrically averaged $\theta$s in two iterations within one time step, 
and the second one uses arithmetically averaged $\theta$s for two consecutive time steps while $\theta$ remains the same in each time step. 
The stability and second order accuracy of our methods are verified using stability and truncation error analysis and are demonstrated by numerical examples on linear and semi-linear hyperbolic partial differential equations and Burgers' equation.
\end{abstract}


\maketitle


\section{Introduction}
\label{sec:intro}
The iterated Crank-Nicolson (ICN) method is a popular and successful numerical method in numerical relativity for solving partial differential equations \cite{Teu00,Lei06}. 
The ICN method is the explicit version of 
the Crank-Nicolson (CN) method, which is a very famous implicit finite difference method for solving partial differential equations \cite{CN}.
The ICN method transforms the implicit CN method into an explicit algorithm through a sequence of iterations. 
It has been suggested by Teukolsky \cite{Teu00} that one should carry out exactly two iterations and no more, 
since the accuracy is not affected by doing more iterations. In this paper we will consider the ICN method with two iterations.
    
The stability of the ICN method can be improved by introducing a variable $\theta$, 
where $\theta$ is the weight when averaging the predicted and corrected values. 
The resulting ICN method is referred as the $\theta$-ICN method, which was introduced by Leiler and Rezzolla in 2006 \cite{Lei06}.
The original ICN method is the special case of $\theta$-ICN when $\theta = 0.5$. 
In numerical relativity simulations, $\theta$ is usually chosen to be larger than 0.5 to obtain better stability.   
For instance, $\theta$ is chosen to be 0.51 in \cite{Yiou09}, and in \cite{Duez03, Duez04} the authors find that $\theta \neq 0.5$ yields an improved stability.  
A major drawback of using $\theta \neq 0.5$ in $\theta$-ICN method is that the accuracy is reduced to first order. Only when $\theta$ = 0.5, the $\theta$-ICN method is second order accurate \cite{Lei06}. 
   
In this paper, we propose two modified $\theta$-ICN algorithms that have second order of convergence when $\theta \neq 0.5$. 
The first one is based on the geometrically averaged weights of two consecutive iterations of the $\theta$-ICN method and this method is referred as the geometric averaging (GA) $\theta$-ICN algorithm.
The second one is based on using two arithmetically averaged weights for two consecutive time steps and this method is referred as the arithmetic averaging (AA) $\theta$-ICN algorithm.
The paper is organized as follows: in section \ref{sec:review}, we review the original ICN and $\theta$-ICN methods. 
In section \ref{sec:method}, we discuss the proposed modified $\theta$-ICN algorithms with improved accuracy. 
Numerical examples on linear hyperbolic PDE, semi-linear hyperbolic PDE, and Burgers' equation are presented in section \ref{sec:examples}.

\section{The ICN and $\theta$-ICN Methods}
\label{sec:review}

Consider the linear hyperbolic PDE
\begin{equation}
 u_{t}+au_{x}=0 , \label{ICND}
\end{equation} 
where $a$ is a constant.
The ICN method solves the implicit Crank-Nicolson update equation by iteration and turn it into an explicit scheme. 
First, the centered difference scheme is used to discretize the equation (\ref{ICND})
\begin{equation}
\dfrac{\tilde{u}^{n+1}_{j}-u^{n}_{j}}{\Delta t}+a \dfrac{u^{n}_{j+1}-u^{n}_{j-1}}{2\Delta x}=0.
\end{equation}  
Solving for $\tilde{u}^{n+1}_{j}$, we get 
\begin{equation}
\tilde{u}^{n+1}_{j}=u^n_j- a \dfrac{\Delta t}{2 \Delta x}(u^{n}_{j+1}-u^{n}_{j-1}).\label{ICN1}
\end{equation} 
Next, an intermediate variable $\bar{u}^{n+1/2}_{j}$ is defined as 
\begin{equation}
\bar{u}^{n+1/2}_{j}=\dfrac{1}{2}(\tilde{u}^{n+1}_{j}+u^{n}_{j}). \label{ICN2}
\end{equation}
Apply the centered difference scheme using intermediate value $\bar{u}^{n+1/2}_{j}$, we get
\begin{equation}
\dfrac{u^{n+1}_{j}-u^{n}_{j}}{\Delta t} + a \dfrac{\bar{u}^{n+1/2}_{j+1}-\bar{u}^{n+1/2}_{j-1}}{2 \Delta x}=0,
\end{equation}
which can be written as
\begin{equation}
u^{n+1}_{j}= u^n_j- a \dfrac{\Delta t}{2 \Delta x}(\bar{u}^{n+1/2}_{j+1}-\bar{u}^{n+1/2}_{j-1}).
\label{ICNe}
\end{equation}
Equations (\ref{ICN1}), (\ref{ICN2}), and (\ref{ICNe}) are the update equations of the ICN method with one iteration. For more iterations, one needs to repeat equations (\ref{ICN1}) and (\ref{ICN2}) .

For the $\theta$-ICN method \cite{Lei06}, we no longer weight $\tilde{u}^{n+1}_{j}$ and $u^{n}_{j}$ equally in equation (\ref{ICN2}). Instead, we define 
\begin{equation}
\bar{u}^{n+1/2}_{j}=\theta \tilde{u}^{n+1}_{j}+(1-\theta)u^{n}_{j},
\end{equation}
where $0 \leq \theta \leq 1$. 

Let $R = a\frac{\Delta t}{2\Delta x}$, solving equation (\ref{ICND}) using $\theta$-ICN method with two iterations, we obtain the following update equations
\begin{align}
^{(1)}\tilde{u}^{n+1}_{j} &= u^{n}_{j}-R(u^{n}_{j+1}-u^{n}_{j-1}) , \\
^{(1)}\bar{u}^{n+1/2}_{j}&=\theta ^{(1)}\tilde{u}^{n+1}_{j}+(1-\theta)u^{n}_{j},\\
^{(2)}\tilde{u}^{n+1}_{j}& = u^n_j - R(^{(1)}\bar{u}^{n+1/2}_{j+1} - ^{(1)}\bar{u}^{n+1/2}_{j-1}), \\
^{(2)}\bar{u}^{n+1/2}_{j}&=\theta ^{(2)}\tilde{u}^{n+1}_{j}+(1-\theta)u^{n}_{j},\\
u^{n+1}_j &= u^n_j - R(^{(2)}\bar{u}^{n+1/2}_{j+1} - ^{(2)}\bar{u}^{n+1/2}_{j-1}).
\end{align}

Leiler and Rezzolla also suggested to swap the weights $\theta$s for $\theta$-ICN \cite{Lei06}. For swapped $\theta$-ICN method, the $\bar{u}^{n+1/2}_{j}$ terms in two consecutive iterations are calculated using swapped weights
\begin{align}
 ^{(1)}\bar{u}^{n+1/2}_{j}=\theta^{(1)}\tilde{u}^{n+1}_{j}+(1-\theta)u^{n}_{j}, \\
 ^{(2)}\bar{u}^{n+1/2}_{j}=(1-\theta)^{(2)}\tilde{u}^{n+1}_{j}+\theta u^{n}_{j}.
\end{align}

\section{Modified $\theta$-ICN Methods}
\label{sec:method}

The $\theta$-ICN method is only first order accurate when $\theta \neq 0.5$. 
In this section, we propose two different ways to modify $\theta$ to achieve second order accuracy when $\theta \neq 0.5$.

\subsection{The Geometric Averaging $\theta$-ICN method}
The first idea is to take the geometric mean of two $\theta$s in two consecutive iterations to be $\frac{1}{2}$. We name this method Geometric Averaging (GA) $\theta$-ICN method. Taking equation (\ref{ICND}) as an example and letting $R=a \frac{\Delta t}{2\Delta x}$, the GA $\theta$-ICN method consists of the following five steps. \\
Step 1. Calculate $\tilde{u}^{n+1}_{j}$
\begin{equation}
\tilde{u}^{n+1}_{j}=u^{n}_{j}-R(u^{n}_{j+1}-u^{n}_{j-1}). 
\end{equation}
Step 2. Average $\tilde{u}^{n+1}_{j}$ and $u^n_j$ using weight $\theta_1$ to obtain $\bar{u}^{n+\theta_1}_{j}$
\begin{equation}
\bar{u}^{n+\theta_1}_{j}=\theta_{1}\tilde{u}^{n+1}_{j}+(1-\theta_{1})u^{n}_{j}.
\end{equation}
Step 3. Calculate $\tilde{u}^{n+2\theta_1}_{j}$ using $\bar{u}^{n+\theta_1}_j$  
\begin{equation}
\tilde{u}^{n+2\theta_1}_{j}=u^{n}_{j}-2\theta_{1}R(\bar{u}^{n+\theta_1}_{j+1}-\bar{u}^{n+\theta_1}_{j-1}) . \label{modify}
\end{equation}
Step 4. Average again using weight $\theta_2$ to obtain $\bar{u}^{n+1/2}_{j}$
\begin{equation}
\bar{u}^{n+1/2}_{j}=\theta_{2} \tilde{u}^{n+2\theta_1}_{j}+(1-\theta_{2})u^{n}_{j}.
\label{GA2}
\end{equation}
Step 5. Compute $u^{n+1}_j$ using $\bar{u}^{n+1/2}_{j}$
\begin{equation}
u^{n+1}_{j}=u^{n}_{j}-R( \bar{u}^{n+1/2}_{j+1}- \bar{u}^{n+1/2}_{j-1}).
\end{equation}
We define $\theta_{1}$ and $\theta_{2}$ to be positive real numbers and their geometric mean to be $\frac{1}{2}$, that is $\sqrt{\theta_1 \theta_2}= \frac{1}{2}$, or equivalently $\theta_1\theta_2=\frac{1}{4}$.

\begin{figure}[ht]
 \centering
 \includegraphics[width=0.8\textwidth]{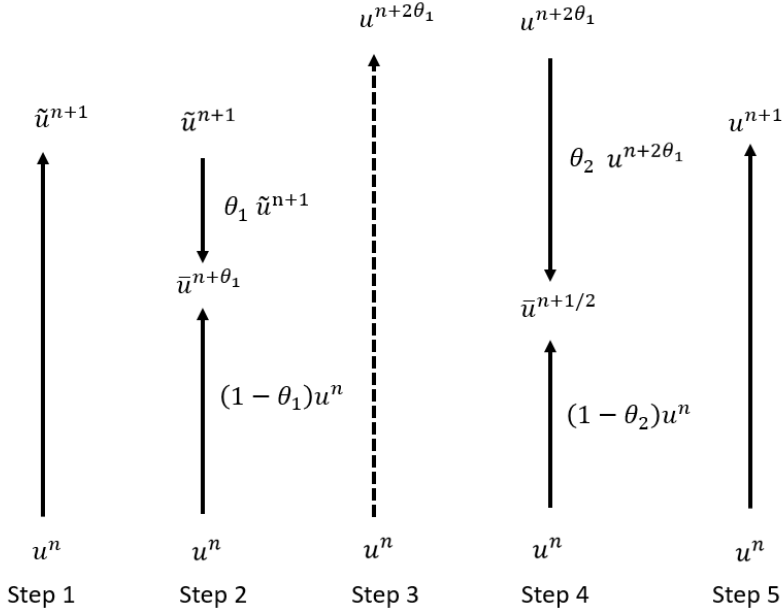}
 \caption{Five steps for the GA $\theta$-ICN method.}
 \label{fig:GAsteps}
\end{figure}

Figure \ref{fig:GAsteps} shows the steps of the GA $\theta$-ICN method. 
In step 2, the averaged value $\bar{u}^{n+\theta_1}_{j}$ is located at time level $n+\theta_1$ (no longer at $n+\frac{1}{2}$ if $\theta_1 \ne 0.5$).
In step 3, the solution is updated from $n$ to $n+2\theta_1$, in order to maintain the centered difference in time. 
Note that equation (\ref{modify}) is obtained from
\begin{equation}
\dfrac{ \tilde{u}^{n+2 \theta_1}_j - u^n_j}{2 \theta_1 \Delta t} + a\dfrac{\bar{u}^{n+ \theta_1}_{j+1}-u^{n+\theta_1}_{j-1}}{2 \Delta x}=0,
\end{equation}
where both $u_t$ and $u_x$ are approximated using centered difference.
In step 4, in order to compute the averaged solution at $n+\frac{1}{2}$ as in equation (\ref{GA2}), we have
\begin{align}
n + \frac{1}{2} = \theta_2 (n+2\theta_1)+(1-\theta_2) n = n+2\theta_1 \theta_2. \label{eq:AAtheta}
\end{align}
Solving equation (\ref{eq:AAtheta}), we get $\theta_1\theta_2=1/4$.
Note that when $\theta_1 = 0.5$ and $\theta_2=1/(4\theta_1)=0.5$, the GA $\theta$-ICN method becomes the standard ICN method ($\theta$-ICN with $\theta=0.5$).

Combining steps 1 to 5, we have
\begin{align}
u^{n+1}_{j}=&-2\theta_{1}^{2}\theta_{2}R^{3}u^{n}_{j+3}+2\theta_{1}\theta_{2}R^{2}u^{n}_{j+2}+(6\theta_{1}^{2}\theta_{2}R^{3}-R) \nonumber \\
&u^{n}_{j+1} +(1-4\theta_{1}\theta_{2}R^{2})u^{n}_{j} + (R-6\theta_{1}^{2}\theta_{2}R^{3})u^{n}_{j-1} \nonumber \\
&+2\theta_{1}\theta_{2}R^{2}u^{n}_{j-2}+2\theta_{1}^{2}\theta_{2}R^{3}u^{n}_{j-3} . \label{theta1u}
\end{align}
Apply the von Neumann stability analysis \cite{Teu00}, and let 
\begin{equation}
u^n_j = \xi^n e^{ikj\Delta x},
\end{equation}
we get the amplification factor
\begin{equation}
g(\xi)=1-2\beta i-2\beta^{2}+4\theta_{1}\beta^{3} i,
\end{equation}
where $\beta=R\sin(k\Delta x)$.

 \begin{figure}[h]
 \centering
  \includegraphics[width=\textwidth]{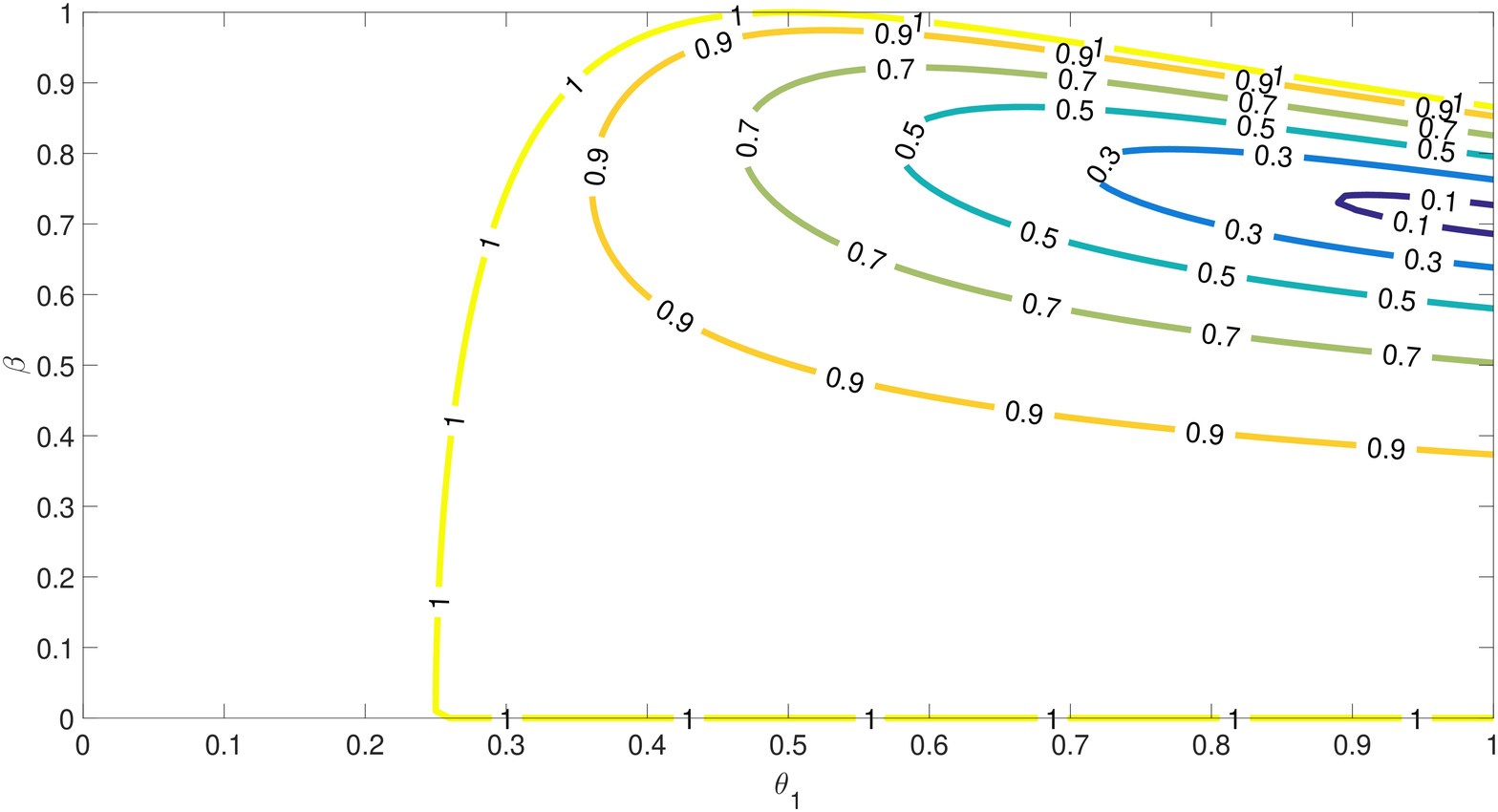}
 \caption{(Color online) Stability Region for the GA $\theta$-ICN method.}
 \label{fig:stGA}
 \end{figure}
Figure \ref{fig:stGA} shows the stability region in the $(\theta_1,\beta)$ plane for the GA $\theta$-ICN method. We see that $\theta_1$ is no longer required to be greater than $0.5$ for the method to be stable.

In the following, we calculate the truncation error of the GA $\theta$-ICN method.
Let
\begin{align}
&\delta^{1}u^{n}_{j}=u^{n}_{j+1}-u^{n}_{j-1}, \\
&\delta^{2}u^{n}_{j}=u^{n}_{j+2}-2u^{n}_{j}+u^{n}_{j-2}, \\
&\delta^{3}u^{n}_{j}=u^{n}_{j+3}-3u^{n}_{j+1}+3u^{n}_{j-1}-u^{n}_{j-3}. 
\end{align}
Equation (\ref{theta1u}) can be rewritten as:
\begin{equation}
u^{n+1}_{j}=u^{n}_{j}-R\delta^{1}u^{n}_{j}+2\theta_{1}\theta_{2}R^{2}\delta^{2}u^{n}_{j}-2\theta_{1}^{2}\theta_{2}R^{3}\delta^{3}u^{n}_{j}. \label{gaunj}
\end{equation}
Use the Taylor expansions, $\delta^{1}u^{n}_{j}$, $\delta^{2}u^{n}_{j}$, and $\delta^{3}u^{n}_{j}$ become
\begin{align}
&\delta^{1}u^{n}_{j}=2 \Delta x (u_x)^n_j+O(\Delta x^3), \label{de1} \\
&\delta^{2}u^{n}_{j}=4 \Delta x^2 (u_{xx})^n_j+O(\Delta x^4),  \label{de2}\\
&\delta^{3}u^{n}_{j}=O(\Delta x^3). \label{de3}
\end{align}
Substituting equations (\ref{de1}), (\ref{de2}) and (\ref{de3}) into equation (\ref{gaunj}), and apply  
$R = a \frac{\Delta t}{2 \Delta x}$, we obtain
\begin{equation}
u^{n+1}_j = u^n_j - a \Delta t(u_x)^n_j + 2 \theta_1 \theta_2 a^2 \Delta t^2 (u_{xx})^n_j + O(\Delta x^3). \label{tt}
\end{equation}
The Taylor expansion of $u^{n+1}_j$ gives
\begin{equation}
u^{n+1}_j = u^n_j + \Delta t (u_t)^n_j  + \frac{\Delta t^2}{2} (u_{tt})^n_j + O(\Delta t^3).
\end{equation}
The resulting truncation error is
\begin{equation}
e_{\tau}=\dfrac{\Delta t}{2}u_{tt}  - 2 \theta_1 \theta_2 a^2 \Delta t u_{xx} +O(\Delta t^2) + O(\Delta x^2). \label{cuo}
\end{equation}
The original differential equation (\ref{ICND}) implies that $u_{tt} = a^2 u_{xx}$, so 
\begin{equation}
e_{\tau} = 2 (\dfrac{1}{4} - \theta_1 \theta_2) a^2 u_{xx} \Delta t + O(\Delta t^2) + O(\Delta x^2).
\label{GAtruncation}
\end{equation}
Since $\theta_1 \theta_2 = \frac{1}{4}$ for the GA $\theta$-ICN method, the first term on the right hand side of the equation (\ref{GAtruncation}) vanishes and the truncation error becomes $e_{\tau} = O(\Delta t^2) + O(\Delta x^2)$. This proves that the GA $\theta$-ICN method is second order accurate in both time and space.

\subsection{The Arithmetic Averaging $\theta$-ICN method}
The second idea is to modify $\theta$ in different time steps.
We define the $\theta$ in odd time steps to be $\theta_{o}$, the $\theta$ in even time steps to be $\theta_{e}$, and we require $\theta_{o}+\theta_{e}=1$. We call this method the Arithmetic Averaging (AA) $\theta$-ICN method, because the arithmetic mean of $\theta_o$ and $\theta_e$ is $\frac{1}{2}$. 

Consider the linear hyperbolic equation (\ref{ICND}) and let $R=a\frac{\Delta t}{2 \Delta x}$.
For the odd time step, we have
\begin{align}
 ^{(1)}\tilde{u}^{n+1}_{j}		&=u^{n}_{j}-R(u^{n}_{j+1}-u^{n}_{j-1}), \label{AAodd1}\\
 ^{(1)}\bar{u}^{n+1/2}_{j} 	&=\theta_{o} ^{(1)}\tilde{u}^{n+1}_{j}+(1-\theta_{o})u^{n}_{j} , \label{AAodd2} \\
 ^{(2)}\tilde{u}^{n+1}_{j}  	&=u^{n}_{j}-R(^{(1)}\bar{u}^{n+1/2}_{j+1}-^{(1)}\bar{u}^{n+1/2}_{j-1}) ,\label{AAodd3} \\
 ^{(2)}\bar{u}^{n+1/2}_{j} 	&=\theta_{o} ^{(2)}\tilde{u}^{n+1}_{j}+(1-\theta_{o})u^{n}_{j} , \label{AAodd4} \\
 u^{n+1}_{j} 				&=u^{n}_{j}-R(^{(2)}\bar{u}^{n+1/2}_{j+1}-^{(2)}\bar{u}^{n+1/2}_{j-1}). \label{AAodd5}
\end{align}
For the next time step (which is an even time step)
\begin{align}
 ^{(1)}\tilde{u}^{n+2}_{j} 	&=u^{n+1}_{j}-R(u^{n+1}_{j+1}-u^{n+1}_{j-1}), \\
 ^{(1)}\bar{u}^{n+3/2}_{j} &=\theta_{e} ^{(1)}\tilde{u}^{n+2}_{j}+(1-\theta_{e})u^{n+1}_{j}, \\
 ^{(2)}\tilde{u}^{n+2}_{j} 	&=u^{n+1}_{j}-R(^{(1)}\bar{u}^{n+3/2}_{j+1}-      ^{(1)}\bar{u}^{n+3/2}_{j-1}) ,\\
 ^{(2)}\bar{u}^{n+3/2}_{j}	&=\theta_{e} ^{(2)}\tilde{u}^{n+2}_{j}+(1-\theta_{e})u^{n+2}_{j}, \\
 u^{n+2}_{j} 			&=u^{n+1}_{j}-R(^{(2)}\bar{u}^{n+3/2}_{j+1}-^{(2)}\bar{u}^{n+3/2}_{j-1}) .
\end{align}
Note that when $\theta_o = \theta_e = 0.5$, the AA $\theta$-ICN method becomes the standard ICN method ($\theta$-ICN with $\theta=0.5$).

Substituting equations (\ref{AAodd1}) - (\ref{AAodd4}) into equation (\ref{AAodd5}), we obtain
\begin{align*}
u^{n+1}_j=&-\theta_{o}^{2}R^{3}u^{n}_{j+3}+\theta_{o}R^{2}u^{n}_{j+2}+(3\theta_{o}^{2}R^{3}-R)u^{n}_{j+1} \\
&+(1-2\theta_{o}R^{2})u^{n}_{j}+(R-3\theta_{o}^{2}R^{3})u^{n}_{j-1} \\
&+\theta_{o}R^{2}u^{n}_{j-2}+\theta_{o}^{2}R^{3}u^{n}_{j-3}. 
\end{align*}
We can find the amplification factor $g_{o}(\xi)$ for the odd time step 
\begin{equation}
g_{o}(\xi)=1-2\beta i-4\theta_{o} \beta^{2}+8\theta_{o}^{2} \beta^{3} i.
\end{equation}
Similarly, the amplification factor $g_{e}(\xi)$ for the even time step is 
\begin{equation}
g_{e}(\xi)=1-2\beta i-4\theta_{e} \beta^{2}+8\theta_{e}^{2} \beta^{3} i .
\end{equation}

 \begin{figure}[h]
 \centering
  \includegraphics[width=\textwidth]{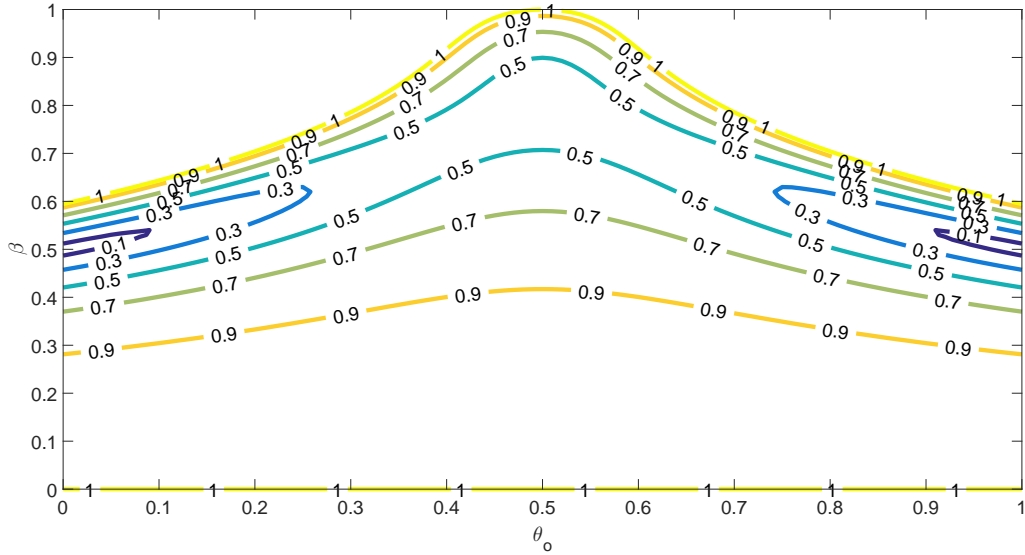}
 \caption{(Color online) Stability Region for the AA $\theta$-ICN method.}
 \label{fig:stAA}
 \end{figure}

The product of two amplification factors $g_{o}(\xi)$ and $g_{e}(\xi)$ is the amplification factor from time steps $n$ to $n+2$. 
Figure \ref{fig:stAA} shows the stability region in $(\theta_o,\beta)$ plane. 
We see that the stability region is symmetric with respect to $\theta_o=0.5$. 
In comparison to the GA $\theta$-ICN method, the AA $\theta$-ICN method has slightly larger stability region, but slightly more damping. 
For example, when $\theta_o=0.4$ and $\beta = 0.6$, the AA $\theta$-ICN method has an amplification factor between $0.5$ and $0.7$, but the GA $\theta$-ICN method has a larger amplification factor around $0.9$.

Following a similar procedure as the GA $\theta$-ICN method, we find that
\begin{equation}
u^{n+1}_j = u^n_j - a\Delta t (u_x)^n_j + a^2 \theta_o \Delta t^2 (u_{xx})^n_j + O(\Delta x^3), \label{un1}
\end{equation}
and
\begin{equation}
u^{n+2}_j = u^{n+1}_j - a \Delta t (u_x)^{n+1}_j + a^2 \theta_e \Delta t^2 (u_{xx})^{n+1}_j + O(\Delta x^3). \label{un2}
\end{equation}
Substitute equation (\ref{un1}) into equation (\ref{un2}), we obtain
\begin{align}
u^{n+2}_j =& u^n_j - a \Delta t (u_x)^n_j + a^2 \theta_o \Delta t^2 (u_{xx})^n_j - a \Delta t (u_x)^{n+1}_j \nonumber \\ 
&+ a^2 \theta_e \Delta t^2 (u_{xx})^{n+1}_j + O(\Delta x^3).  \label{un1un2}
\end{align}
The Taylor expansions of $(u_x)^{n+1}_j$ and $(u_{xx})_j^{n+1}$ give
\begin{align}
&(u_x)^{n+1}_j = (u_x)^n_j + \Delta t (u_{xt})^n_j + O(\Delta t^2), \label{taylorx} \\
&(u_{xx})^{n+1}_j = (u_{xx})^n_j + O(\Delta t). \label{taylorxx}
\end{align}
Substitute equations (\ref{taylorx}) and (\ref{taylorxx}) into equation (\ref{un1un2}), we get
\begin{align}
u^{n+2}_j =& u^n_j - 2 a \Delta t (u_x)^n_j + a^2 \theta_o \Delta t^2 (u_{xx})^n_j - a \Delta t^2 (u_{xt})^n_j  \nonumber \\
&+ a^2 \theta_e \Delta t^2 (u_{xx})^n_j + O(\Delta t^3) +O(\Delta x^3).
\end{align}
The original differential equation (\ref{ICND}) implies that $u_{xt} = -a u_{xx}$, so we have
\begin{align}
u^{n+2}_j = & u^n_j -2 a \Delta t (u_x)^n_j + a^2 \Delta t^2 (1+\theta_o +\theta_e)  (u_{xx})^n_j  \nonumber \\
&+ O(\Delta x^3) +O(\Delta t^3). \label{tt}
\end{align}
The Taylor expansion of $u^{n+2}_j$ gives
\begin{equation}
u^{n+2}_j = u^n_j +2 \Delta t u_t +\dfrac{(2\Delta t)^2}{2} u_{tt} + O(\Delta t^3). \label{xx}
\end{equation}
Substituting equation (\ref{xx}) into equation (\ref{tt}), we get
\begin{equation}
u_t+ a u_x = -\Delta t u_{tt} + \dfrac{a^2 \Delta t}{2}(1+\theta_o+\theta_e)u_{xx}+O(\Delta t^2)+ O(\Delta x^2)  .\label{11}
\end{equation}
The original differential equation (\ref{ICND}) implies that $u_{tt} = a^2 u_{xx}$, so we obtain the truncation error
\begin{equation}
e_{\tau} = (\dfrac{1+\theta_o+\theta_e}{2}-1) a^2 u_{xx} \Delta t +O(\Delta t^2) +O(\Delta x^2).
\label{AAtruncation}
\end{equation}
Since the AA $\theta$-ICN method requires $\theta_o+\theta_e=1$, the first term on the right hand side of the equation (\ref{AAtruncation}) vanishes. Therefore, the truncation error becomes $e_{\tau} = O(\Delta t^2) +O(\Delta x^2)$, which indicates that the AA $\theta$-ICN method with two iterations is second order accurate in both time and space.

\section{Numerical Examples}
\label{sec:examples}

\subsection{Linear hyperbolic PDE}
In the first numerical example, we consider the following linear hyperbolic equation initial value problem with periodic boundary condition 
\begin{align}
&u_{t}+u_{x}=0, x\in[0,1], t\in[0,1], \\
&u(0,x)=\sin^{2}(\pi x), u(x+1,t)=u(x,t).
\end{align}
The exact solution is $u(x,t)=\sin^{2}((x-t)\pi)$. The Courant-Friedrichs-Lewy (CFL) condition is chosen to be $0.5$ and the numbers of grid points in space are chosen to be $N=100,200,400,800,$ and $1600$. For the GA $\theta$-ICN method, we chose $\theta_{1}=0.6$ and $\theta_{2}=\frac{1}{4\theta_1} \approx 0.416$. For the AA $\theta$-ICN method, we use $\theta_o =0.6$ and $\theta_e = 0.4$. We compare our new methods with the standard ICN method, the $\theta$-ICN method with $\theta = 0.6$, and the swapped $\theta$-ICN method. 

\begin{figure}[h]
 \centering
 \includegraphics[width=\textwidth]{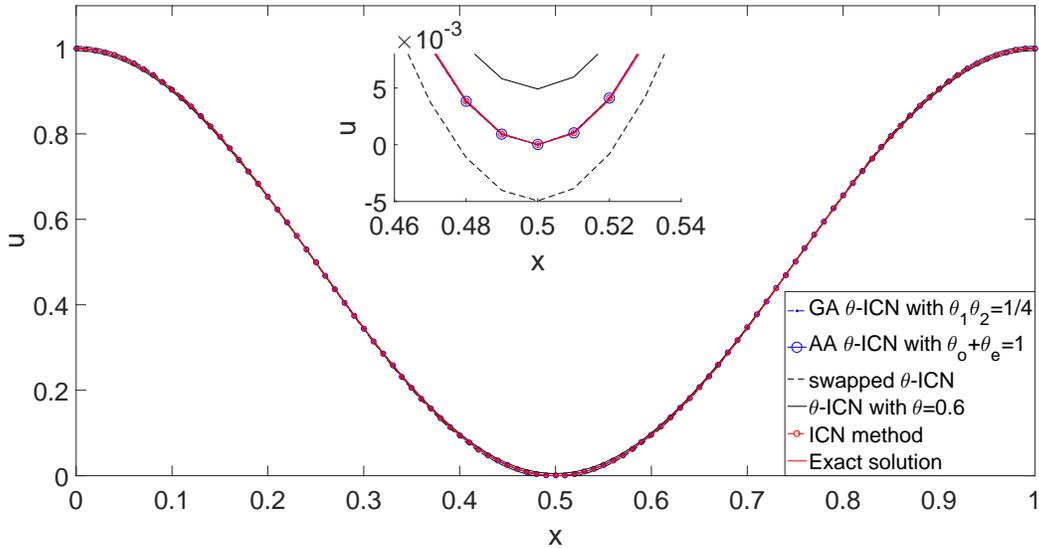}
 \caption{(Color online) Comparison between the GA $\theta$-ICN method, the AA $\theta$-ICN method, the swapped $\theta$-ICN method, the $\theta$-ICN method ($\theta = 0.6$), the ICN method, and the exact solution for linear hyperbolic PDE.}
 \label{fig:linear}
\end{figure}

Figure \ref{fig:linear} shows the numerical results and the exact solution. We see that the GA $\theta$-ICN method, the AA $\theta$-ICN method, and the ICN method are more accurate than the swapped $\theta$-ICN method and the $\theta$-ICN method with $\theta = 0.6$.
We calculate the $L_1$, $L_2$, and $L_\infty$ norms of each method and the results are shown in tables \ref{table:1}, \ref{table:2}, and \ref{table:3}, respectively. 
From these tables, we see that the swapped $\theta$-ICN method and the $\theta$-ICN method with $\theta=0.6$ are only first order accurate, while the ICN method, the GA and AA $\theta$-ICN methods are second order accurate.

\begin{table}[ht]
\caption{Comparison of $L_1$ Norm for Linear Hyperbolic PDE.}
\begin{center}
\small
\begin{tabular}{|c|c|c|c|c|c|c|c|c|c|c|}
\hline
&\multicolumn{2}{c|}{ICN}&\multicolumn{2}{c|}{$\theta$-ICN}&\multicolumn{2}{c|}{Swapped $\theta$-ICN}&\multicolumn{2}{c|}{GA $\theta$-ICN}&\multicolumn{2}{c|}{AA $\theta$-ICN}\\
\cline{2-11}
&\multicolumn{2}{c|}{$\theta=0.5$}&\multicolumn{2}{c|}{$\theta=0.6$}&\multicolumn{2}{c|}{$\theta_1+\theta_2=1$}&\multicolumn{2}{c|}{$\theta_{1}\theta_{2}=1/4$}&\multicolumn{2}{c|}{$\theta_{o}+\theta_e=1$}\\
\hline
N&$L_1$&order&$L_1$&order&$L_1$&order&$L_1$&order&$L_1$&order\\
\hline
200&1.8E-4& &1.6E-3& &1.6E-3& &2.0E-4& &1.9E-4&\\
\hline
400&4.6E-5&2.0&7.9E-4&1.0&7.9E-4&1.0&4.9E-5&2.0&4.7E-5&2.0\\
\hline
800&1.2E-5&2.0&3.9E-4&1.0&3.9E-4&1.0&1.2E-5&2.0&1.2E-5&2.0\\
\hline
1600&2.9E-6&2.0&2.0E-4&1.0&2.0E-4&1.0&3.1E-6&2.0&2.9E-6&2.0\\
\hline
\end{tabular}
\end{center}
\label{table:1}
\end{table}

\begin{table}[ht]
\caption{Comparison of $L_2$ Norm for Linear Hyperbolic PDE.}
\begin{center}
\small
\begin{tabular}{|c|c|c|c|c|c|c|c|c|c|c|}
\hline
&\multicolumn{2}{c|}{ICN}&\multicolumn{2}{c|}{$\theta$-ICN}&\multicolumn{2}{c|}{Swapped $\theta$-ICN}&\multicolumn{2}{c|}{GA $\theta$-ICN}&\multicolumn{2}{c|}{AA $\theta$-ICN}\\
\cline{2-11}
&\multicolumn{2}{c|}{$\theta=0.5$}&\multicolumn{2}{c|}{$\theta=0.6$}&\multicolumn{2}{c|}{$\theta_1+\theta_2=1$}&\multicolumn{2}{c|}{$\theta_{1} \theta_{2}=1/4$}&\multicolumn{2}{c|}{$\theta_{o}+\theta_e=1$}\\
\hline
N&$L_2$&order&$L_2$&order&$L_2$&order&$L_2$&order&$L_2$&order\\
\hline
200&1.5E-5& &1.3E-4& &1.3E-4& &1.5E-5& &1.5E-5&\\
\hline
400&2.6E-6&2.5&4.4E-5&1.5&4.4E-5&1.5&2.7E-6&2.5&2.6E-6&2.5\\
\hline
800&4.5E-7&2.5&1.5E-5&1.5&1.6E-5&1.5&4.8E-7&2.5&4.6E-7&2.5\\
\hline
1600&8.0E-8&2.5&5.5E-6&1.5&5.5E-6&1.5&8.6E-8&2.5&8.1E-8&2.5\\
\hline
\end{tabular}
\end{center}
\label{table:2}
\end{table}

\begin{table}[ht]
\caption{Comparison of $L_{\infty}$ Norm for Linear Hyperbolic PDE.}
\begin{center}
\small
\begin{tabular}{|c|c|c|c|c|c|c|c|c|c|c|}
\hline
&\multicolumn{2}{c|}{ICN}&\multicolumn{2}{c|}{$\theta$-ICN}&\multicolumn{2}{c|}{Swapped $\theta$-ICN}&\multicolumn{2}{c|}{GA $\theta$-ICN}&\multicolumn{2}{c|}{AA $\theta$-ICN}\\
\cline{2-11}
&\multicolumn{2}{c|}{$\theta=0.5$}&\multicolumn{2}{c|}{$\theta=0.6$}&\multicolumn{2}{c|}{$\theta_1+\theta_2=1$}&\multicolumn{2}{c|}{$\theta_{1} \theta_{2}=1/4$}&\multicolumn{2}{c|}{$\theta_{o}+\theta_e=1$}\\
\hline
N&$L_{\infty}$&order&$L_{\infty}$&order&$L_{\infty}$&order&$L_{\infty}$&order&$L_{\infty}$&order\\
\hline
200&2.9E-4& &2.5E-3& &2.5E-3& &3.1E-4& &3.0E-4&\\
\hline
400&7.3E-5&2.0&1.2E-3&1.1&1.2E-3&1.1&7.8E-5&2.0&7.4E-5&2.0\\
\hline
800&1.8E-5&2.0&6.2E-4&1.0&6.2E-4&1.0&2.0E-5&2.0&1.8E-5&2.0  \\
\hline
1600&4.5E-6&2.0&3.1E-4&1.0&3.2E-4&1.0&4.8E-6&2.0&4.6E-6&2.0\\
\hline
\end{tabular}
\end{center}
\label{table:3}
\end{table}

\subsection{Semi-linear hyperbolic PDE}
In the second example, we consider a semi-linear hyperbolic PDE 
\begin{align}
&u_{t}+u_{x}=-u^{2}, x\in[0,1], t\in[0,1],  \\
&u(0,x)=\sin^{2}(\pi x), u(x+1,t)=u(x,t). \label{semilinear}
\end{align}
The exact solution is given by \cite{Bour10}
\begin{equation}
u(x,t) = \dfrac{\sin^2((x-t)\pi)}{1+t \sin^2((x-t)\pi)}.
\end{equation}
We use the same grid and the CFL condition as in the previous example.  
We solve this PDE by the GA $\theta$-ICN method with $\theta_{1}=0.6$, the AA $\theta$-ICN method with $\theta_{o}=0.6$, the $\theta$-ICN method with $\theta=0.6$, the swapped $\theta$-ICN method, and the ICN method. 

\begin{figure}[h]
 \centering
 \includegraphics[width=\textwidth]{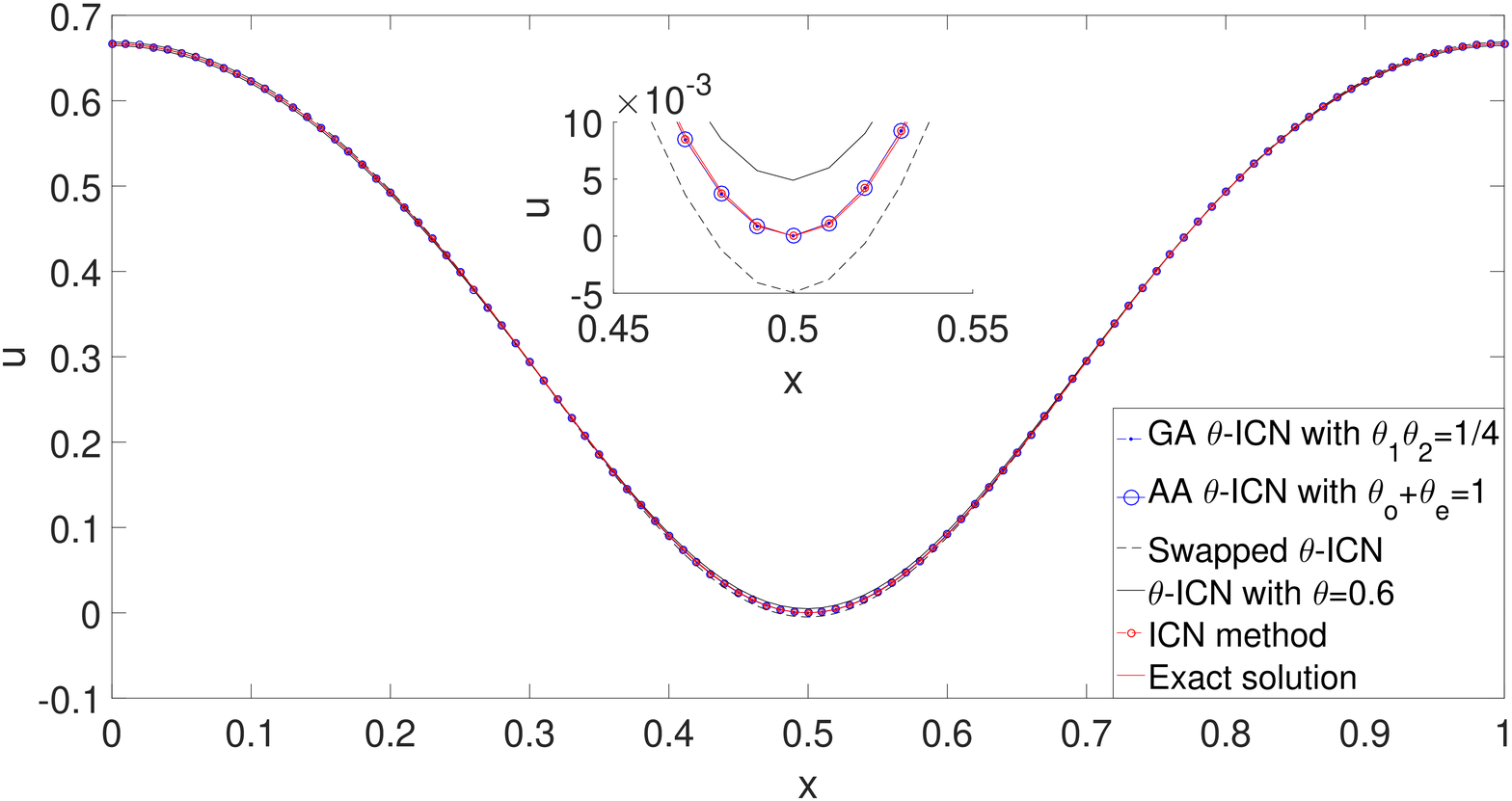}
 \caption{(Color online) Comparison between the GA $\theta$-ICN method, the AA $\theta$-ICN method, the $\theta$-ICN method ($\theta = 0.6$), the swapped $\theta$-ICN method, the ICN method, and the exact solution for Semi-linear Hyperbolic PDE.}
 \label{fig:semilinear}
\end{figure}

Figure \ref{fig:semilinear} shows the results and we see that all methods give correct solutions. 
From the enlarged figure, we can see that the GA and AA $\theta$-ICN methods have similar accuracy as the ICN method, while they are more accurate than the $\theta$-ICN method with $\theta=0.6$ and swapped $\theta$-ICN method. 
We also calculate the $L_1$, $L_2$, and $L_{\infty}$ norms to test the convergence of our proposed methods.
Tables \ref{table:4}, \ref{table:5}, and \ref{table:6} show the $L_1$, $L_2$, and $L_{\infty}$ norm, respectively. From these tables, we see that the GA and AA $\theta$-ICN methods are second order accurate. In comparison, the swapped $\theta$-ICN method and the $\theta$-ICN method with $\theta=0.6$ are only first order accurate.

\begin{table}[ht]
\caption{Comparison of $L_{1}$ Norm for Semi-Linear Hyperbolic PDE.}
\begin{center}
\small
\begin{tabular}{|c|c|c|c|c|c|c|c|c|c|c|}
\hline
&\multicolumn{2}{c|}{ICN}&\multicolumn{2}{c|}{$\theta$-ICN}&\multicolumn{2}{c|}{Swapped $\theta$-ICN}&\multicolumn{2}{c|}{GA $\theta$-ICN}&\multicolumn{2}{c|}{AA $\theta$-ICN}\\
\cline{2-11}
&\multicolumn{2}{c|}{$\theta=0.5$}&\multicolumn{2}{c|}{$\theta=0.6$}&\multicolumn{2}{c|}{$\theta_1+\theta_2=1$}&\multicolumn{2}{c|}{$\theta_{1}\theta_{2}=1/4$}&\multicolumn{2}{c|}{$\theta_{o}+\theta_e=1$}\\
\hline
N&$L_1$&order&$L_1$&order&$L_1$&order&$L_1$&order&$L_1$&order\\
\hline
200&1.3E-4& &1.1E-3& &1.1E-3& &1.4E-4& &1.3E-4&\\
\hline
400&3.3E-5&2.0&5.4E-4&1.0&5.3E-4&1.1&3.5E-5&2.0&3.3E-5&2.0\\
\hline
800&8.1E-6&2.0&2.7E-4&1.0&2.7E-4&1.0&8.7E-6&2.0&8.2E-6&2.0\\
\hline
1600&2.0E-6&2.0&1.4E-4&1.0&1.3E-4&1.0&2.2E-6&2.0&2.1E-6&2.0\\
\hline
\end{tabular}
\end{center}
\label{table:4}
\end{table}

\begin{table}[ht]
\caption{Comparison of $L_2$ Norm for Semi-Linear Hyperbolic PDE.}
\begin{center}
\small
\begin{tabular}{|c|c|c|c|c|c|c|c|c|c|c|}
\hline
&\multicolumn{2}{c|}{ICN}&\multicolumn{2}{c|}{$\theta$-ICN}&\multicolumn{2}{c|}{Swapped $\theta$-ICN}&\multicolumn{2}{c|}{GA $\theta$-ICN}&\multicolumn{2}{c|}{AA $\theta$-ICN}\\
\cline{2-11}
&\multicolumn{2}{c|}{$\theta=0.5$}&\multicolumn{2}{c|}{$\theta=0.6$}&\multicolumn{2}{c|}{$\theta_1+\theta_2=1$}&\multicolumn{2}{c|}{$\theta_{1}\theta_{2}=1/4$}&\multicolumn{2}{c|}{$\theta_{o}+\theta_e=1$}\\
\hline
N&$L_2$&order&$L_2$&order&$L_2$&order&$L_2$&order&$L_2$&order\\
\hline
200&1.1E-5& &9.2E-5& &8.9E-5& &1.2E-5& &1.1E-5&\\
\hline
400&2.0E-6&2.5&3.2E-5&1.5&3.2E-5&1.5&2.1E-6&2.5&2.0E-6&2.5\\
\hline
800&3.5E-7&2.5&1.1E-5&1.5&1.1E-5&1.5&3.7E-7&2.5&3.5E-7&2.5\\
\hline
1600&6.1E-8&2.5&4.0E-6&1.5&4.0E-6&1.5&6.5E-8&2.5&6.2E-8&2.5\\
\hline
\end{tabular}
\end{center}
\label{table:5}
\end{table}

\begin{table}[ht]
\caption{Comparison of $L_{\infty}$ Norm for Semi-Linear Hyperbolic PDE.}
\begin{center}
\small
\begin{tabular}{|c|c|c|c|c|c|c|c|c|c|c|}
\hline
&\multicolumn{2}{c|}{ICN}&\multicolumn{2}{c|}{$\theta$-ICN}&\multicolumn{2}{c|}{Swapped $\theta$-ICN}&\multicolumn{2}{c|}{GA $\theta$-ICN}&\multicolumn{2}{c|}{AA $\theta$-ICN}\\
\cline{2-11}
&\multicolumn{2}{c|}{$\theta=0.5$}&\multicolumn{2}{c|}{$\theta=0.6$}&\multicolumn{2}{c|}{$\theta_1+\theta_2=1$}&\multicolumn{2}{c|}{$\theta_{1}\theta_{2}=1/4$}&\multicolumn{2}{c|}{$\theta_{o}+\theta_e=1$}\\
\hline
N&$L_\infty$&order&$L_\infty$&order&$L_\infty$&order&$L_\infty$&order&$L_\infty$&order\\
\hline
200&2.7E-4& &2.5E-3& &2.5E-3& &2.9E-4& &2.7E-4&\\
\hline
400&6.7E-5&2.0&1.2E-3&1.1&1.2E-3&1.1&7.2E-5&2.0&6.8E-5&2.0\\
\hline
800&1.7E-5&2.0&6.2E-4&1.0&6.2E-4&1.0&1.8E-5&2.0&1.7E-5&2.0\\
\hline
1600&4.2e-6&2.0&3.1E-4&1.0&3.2E-4&1.0&4.5E-6&2.0&4.3E-6&2.0 \\
\hline
\end{tabular}
\end{center}
\label{table:6}
\end{table}

\subsection{Burgers' equation}

In the third example, we consider the Burgers' equation
\begin{align}
&u_{t}+uu_{x} = au_{xx}, x\in[0,1], t\in[0,1], \\
&u(x,0)=\sin^{2}(\pi x), u(x+1,t)=u(x,t),
\end{align}
where $a$ is chosen to be 0.01 in our simulation. 
The term $uu_x$ can be written in conservation form $(F(u))_x$, where $F(u)=\frac{1}{2}u^2$. 
The term $u_{xx}$ can be approximated by the centered difference
\begin{equation}
u_{xx} = \dfrac{u^n_{j+1}-2u^n_j+u^n_{j-1}}{\Delta x^2}.
\end{equation}

\begin{figure}[h]
 \centering
 \includegraphics[width=\textwidth]{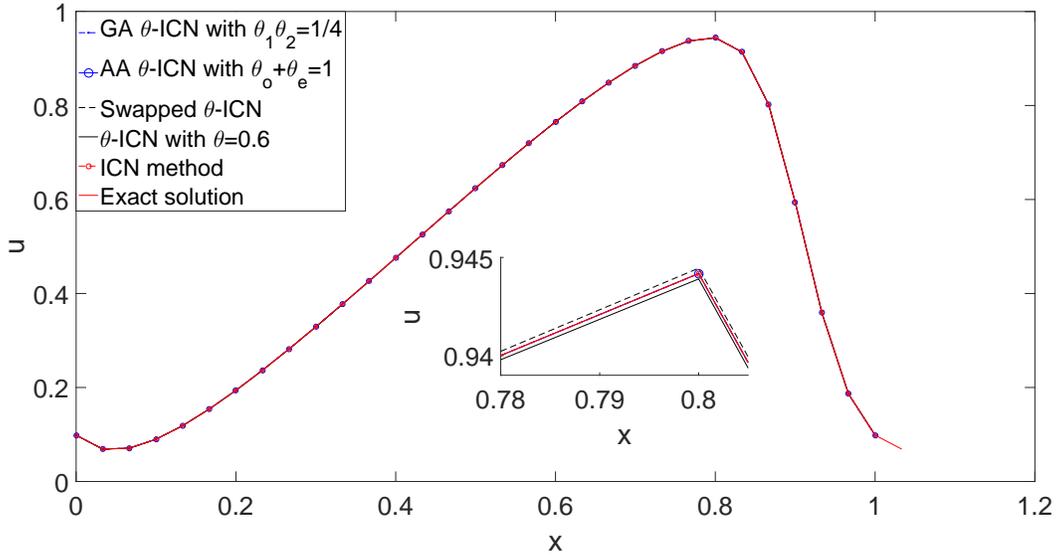}
 \caption{(Color online) Comparison between the GA $\theta$-ICN method, the AA $\theta$-ICN method, the swapped $\theta$-ICN method, the $\theta$-ICN method ($\theta = 0.6$), the ICN method, and the exact solution for Burgers' Equation.}
 \label{fig:burgurs}
\end{figure}

We let the grid size to be $N=30$.
To test the convergence rate in temporal domain, we refine $\Delta t$ and keep $\Delta x$ the same in our simulations, so we let  $\Delta t$ to be $\Delta$, $\Delta/2$, $\Delta/4$, and $\Delta/8$, where $\Delta = 0.5\Delta x^2$.
We use the ICN method with $\Delta t=\Delta/32$ as the exact solution. 
Similar to previous examples, we compare the numerical results of the GA $\theta$-ICN method with $\theta_1=0.6$, the AA $\theta$-ICN method with $\theta_o=0.6$, the $\theta$-ICN method with $\theta=0.6$, the swapped $\theta$-ICN method, and the ICN method. 
Figure \ref{fig:burgurs} shows the results and we see that the GA and AA $\theta$-ICN methods are very close to the ICN method and they are more accurate than the $\theta$-ICN with $\theta=0.6$ and the swapped $\theta$-ICN method.
Tables \ref{table:7}, \ref{table:8} and \ref{table:9} show the numerical results on the $L_1$, $L_2$ and $L_{\infty}$ norms,  respectively. We see that the GA and the AA $\theta$-ICN methods are second order accurate in time, while the swapped $\theta$-ICN method and the $\theta$-ICN method with $\theta$=0.6 are only first order accurate in time. From this example, we see that the GA and AA $\theta$-ICN methods are suitable for solving nonlinear mixed hyperbolic-parabolic equations with improved accuracy.

\begin{table}[ht]
\caption{Comparison of $L_{1}$ Norm in time for Burgers' equation.}
\begin{center}
\small
\begin{tabular}{|c|c|c|c|c|c|c|c|c|c|c|}
\hline
&\multicolumn{2}{c|}{ICN}&\multicolumn{2}{c|}{$\theta$-ICN}&\multicolumn{2}{c|}{Swapped $\theta$-ICN}&\multicolumn{2}{c|}{GA $\theta$-ICN}&\multicolumn{2}{c|}{AA $\theta$-ICN}\\
\cline{2-11}
&\multicolumn{2}{c|}{$\theta=0.5$}&\multicolumn{2}{c|}{$\theta=0.6$}&\multicolumn{2}{c|}{$\theta_1+\theta_2=1$}&\multicolumn{2}{c|}{$\theta_{1}\theta_{2}=1/4$}&\multicolumn{2}{c|}{$\theta_{o}+\theta_e=1$}\\
\hline
$\Delta t$&$L_1$&order&$L_1$&order&$L_1$&order&$L_1$&order&$L_1$&order\\
\hline
$\Delta$&2.9E-7& &7.8E-5& &7.8E-5& &4.7E-7& &3.4E-7&\\
\hline
$\Delta/2$&7.3E-8&2.0&3.9E-5&1.0&3.9E-5&1.0&1.2E-7&2.0&8.5E-8&2.0\\
\hline
$\Delta/4$&1.8E-8&2.0&1.9E-5&1.0&1.9E-5&1.0&2.9E-8&2.0&2.1E-8&2.0\\
\hline
$\Delta/8$&4.3E-9&2.1&9.7E-6&1.0&9.7E-6&1.0&7.1E-9&2.0&5.0E-9&2.1\\
\hline
\end{tabular}
\end{center}
\label{table:7}
\end{table}

\begin{table}[ht]
\caption{Comparison of $L_2$ Norm in time for Burgers' equation.}
\begin{center}
\small
\begin{tabular}{|c|c|c|c|c|c|c|c|c|c|c|}
\hline
&\multicolumn{2}{c|}{ICN}&\multicolumn{2}{c|}{$\theta$-ICN}&\multicolumn{2}{c|}{Swapped $\theta$-ICN}&\multicolumn{2}{c|}{GA $\theta$-ICN}&\multicolumn{2}{c|}{AA $\theta$-ICN}\\
\cline{2-11}
&\multicolumn{2}{c|}{$\theta=0.5$}&\multicolumn{2}{c|}{$\theta=0.6$}&\multicolumn{2}{c|}{$\theta_1+\theta_2=1$}&\multicolumn{2}{c|}{$\theta_{1}\theta_{2}=1/4$}&\multicolumn{2}{c|}{$\theta_{o}+\theta_e=1$}\\
\hline
$\Delta t$&$L_2$&order&$L_2$&order&$L_2$&order&$L_2$&order&$L_2$&order\\
\hline
$\Delta$&9.0E-8& &2.0E-5& &2.0E-5& &1.4E-7& &1.0E-7&\\
\hline
$\Delta/2$&2.3E-8&2.0&1.0E-5&1.0&1.0E-5&1.0&3.6E-8&2.0&2.6E-8&1.9\\
\hline
$\Delta/4$&5.6E-9&2.0&5.0E-6&1.0&5.0E-6&1.0&8.9E-9&2.0&6.3E-9&2.0\\
\hline
$\Delta/8$&1.3E-9&2.1&2.5E-6&1.0&2.5E-6&1.0&2.2E-9&2.0&1.5E-9&2.1\\
\hline
\end{tabular}
\end{center}
\label{table:8}
\end{table}

\begin{table}[ht]
\caption{Comparison of $L_{\infty}$ Norm in time for Burgers' equation.}
\begin{center}
\small
\begin{tabular}{|c|c|c|c|c|c|c|c|c|c|c|}
\hline
&\multicolumn{2}{c|}{ICN}&\multicolumn{2}{c|}{$\theta$-ICN}&\multicolumn{2}{c|}{Swapped $\theta$-ICN}&\multicolumn{2}{c|}{GA $\theta$-ICN}&\multicolumn{2}{c|}{AA $\theta$-ICN}\\
\cline{2-11}
&\multicolumn{2}{c|}{$\theta=0.5$}&\multicolumn{2}{c|}{$\theta=0.6$}&\multicolumn{2}{c|}{$\theta_1+\theta_2=1$}&\multicolumn{2}{c|}{$\theta_{1}\theta_{2}=1/4$}&\multicolumn{2}{c|}{$\theta_{o}+\theta_e=1$}\\
\hline
$\Delta t$&$L_\infty$&order&$L_\infty$&order&$L_\infty$&order&$L_\infty$&order&$L_\infty$&order\\
\hline
$\Delta $&1.7E-6& &3.3E-4& &3.4E-4& &2.7E-6& &1.8E-6&\\
\hline
$\Delta/2$&4.2E-7&2.0&1.7E-4&1.0&1.7E-4&1.0&6.7E-7&2.0&4.6E-7&2.0\\
\hline
$\Delta/4$&1.0E-7&2.1&8.4E-5&1.0&8.4E-5&1.0&1.7E-7&2.0&1.1E-7&2.1\\
\hline
$\Delta/8$&2.5E-8&2.0&4.2E-5&1.0&4.2E-5&1.0&4.0E-8&2.1&2.7E-8&2.0\\
\hline
\end{tabular}
\end{center}
\label{table:9}
\end{table}

\section{Conclusion}
\label{sec:conc}
In this paper, we have proposed two approaches to improve the $\theta$-iterated Crank-Nicolson (ICN) method to second order accuracy when $\theta$ does not equal to 0.5. 
The first approach employs geometrically averaged $\theta$s in two iterations within one time step. 
The second approach uses arithmetically averaged $\theta$s for two consecutive time steps while same $\theta$ is used in every iteration of each time step. 
Stability and truncation error analysis have been carried out to show that our methods are stable and second order accurate.
Numerical examples on linear hyperbolic PDE, semi-linear hyperbolic PDE, and Burgers' equation are presented to verify that the second order accuracy of the proposed new methods.

\begin{acknowledgments}
This work was supported in part by 
the AFOSR Grant FA9550-16-1-0199, 
the US ARO Grant W911NF-11-2-0046, 
and the NSF Grant HRD-1242067.
\end{acknowledgments}

\bibliography{refs}

\end{document}